\documentclass[11pt]{article}
\pdfoutput=1
\usepackage{amsmath,amsthm,amssymb,amsfonts,bbm,bm,math}
\usepackage{algorithmic}
\usepackage[ruled,vlined]{algorithm2e}
\usepackage{accents}
\usepackage{epsfig,epsf,url,enumitem}
\usepackage{etoolbox}
\usepackage{graphicx}
\usepackage{fullpage}
\usepackage{subcaption}

\usepackage{pict2e}

\usepackage{microtype}

\usepackage{booktabs} % for professional tables
\usepackage[dvipsnames]{xcolor}
\usepackage{caption}
\usepackage{colortbl}

\usepackage{mathtools}
\usepackage{times}
\usepackage{latexsym}
\usepackage{multirow}
\usepackage{tabularx}
 \usepackage[square,numbers]{natbib}
\usepackage[colorlinks=true, citecolor=blue]{hyperref} 

\usepackage[most]{tcolorbox}

\newtcolorbox{conversationbox}{sharp corners, colback=gray!10, colframe=gray!60, boxrule=0.1mm, boxsep=0mm, breakable}
\newtcolorbox{humaninput}{colback=blue!5, boxrule=0mm, boxsep=1mm, left=1mm, right=1mm, breakable, before upper=\strut, top=1mm, bottom=1mm}
\newtcolorbox{gptoutput}{colback=green!5, boxrule=0mm, boxsep=1mm, left=1mm, right=1mm, breakable, before upper=\strut, top=1mm, bottom=1mm}
\newtcolorbox{claude3output}{colback=yellow!5, boxrule=0mm, boxsep=1mm, left=1mm, right=1mm, breakable, before upper=\strut, top=1mm, bottom=1mm}
\newtcolorbox{geminioutput}{colback=red!5, boxrule=0mm, boxsep=1mm, left=1mm, right=1mm, breakable, before upper=\strut, top=1mm, bottom=1mm}

\usepackage{fvextra}
\usepackage[frozencache=true, finalizecache=false, cachedir=./minted-cache]{minted} 
\usepackage{float}
\usepackage{alltt}
\usepackage{soul}

\usepackage{url}
\usepackage{MnSymbol}
\usepackage{verbatim}
\usepackage{fontawesome}

\usepackage[T1]{fontenc}
\usepackage[utf8]{inputenc}
\usepackage{varwidth}
\usepackage{listings}

    \usepackage[most]{tcolorbox}
    
    \usepackage{tikz}
    \usetikzlibrary{shapes,calc,positioning}
    \tcbset{
      aibox/.style={
        width=474.18663pt,
        top=10pt,
        colback=white,
        colframe=black,
        colbacktitle=black,
        enhanced,
        center,
        attach boxed title to top left={yshift=-0.1in,xshift=0.15in},
        boxed title style={boxrule=0pt,colframe=white,},
      }
    }
    \newtcolorbox{AIbox}[2][]{aibox,title=#2,#1}
    
    \definecolor{aigold}{RGB}{244,210, 1} 
    \definecolor{aigreen}{RGB}{210,244,211} 
    
    \sethlcolor{aigreen}
    
    \definecolor{aired}{RGB}{255,180,181}

    \newtcbox{\mybox}[1][green]{on line,
    arc=0pt,outer arc=0pt,colback=#1!10!white,colframe=#1!50!black,
    boxsep=0pt,left=0pt,right=0pt,top=0pt,bottom=0pt,
    boxrule=0pt,bottomrule=0pt,toprule=0pt}

    \usepackage[textsize=tiny]{todonotes}
    \usepackage{pifont} % Add this line for \ding commands

    \usepackage{makecell}
    \usepackage{array}
    \newcolumntype{C}[1]{>{\centering\let\newline\\\arraybackslash\hspace{0pt}}m{#1}}
    
    \hypersetup{
      colorlinks   = true, %Colours links instead of ugly boxes
      urlcolor     = blue, %Colour for external hyperlinks
      linkcolor    = blue, %Colour of internal links
      citecolor   = blue %Colour of citations
    }

    \usepackage{longtable}
    
    \usepackage{amsbsy}
    \usepackage{bbm}

    \newcommand{\bsmtx}{\left[ \begin{smallmatrix}} 
    \newcommand{\esmtx}{\end{smallmatrix} \right]} 
    \newcommand{\bmatarray}[1]{\left[\begin{array}{#1}}
    \newcommand{\ematarray}{\end{array}\right]}

    \makeatletter
    \newcommand*{\rom}[1]{\expandafter\@slowromancap\romannumeral #1@}
    \makeatother

    \newcommand\EnumPrefix{}
    
    \newlist{senenum}{enumerate}{10}
    \setlist[senenum]{label=\arabic*.,ref=\EnumPrefix,leftmargin=*}
    
    \AtBeginEnvironment{lemmalist}{\renewcommand\EnumPrefix{\thelemmalist.\arabic*}}

    \numberwithin{equation}{section}

    \newcommand{\remove}[1]{}

\title{\Large\bf Capabilities of Large Language Models in Control Engineering:\\  A Benchmark Study on GPT-4, Claude 3 Opus, and~Gemini~1.0~Ultra}

\author{Darioush Kevian$^1$, Usman Syed$^{1}$, Xingang Guo$^{1}$, Aaron Havens$^{1}$,\\ Geir Dullerud$^{1}$, Peter Seiler$^{2}$,
Lianhui Qin$^{34}$,  and Bin Hu$^{*1}$
\normalsize
}
\date{%
    $^1$University of Illinois Urbana-Champaign\\%
    $^2$University of Michigan\\
    $^3$Allen Institute for AI \\
    $^4$University of California San Diego\\
  $^*$ Corresponding author. E-Mail: binhu7@illinois.edu
}

\begin{document}

\maketitle

\begin{abstract}
In this paper, we explore the capabilities of state-of-the-art large language models (LLMs) such as GPT-4, Claude 3 Opus, and Gemini 1.0 Ultra in solving undergraduate-level control problems. Controls provides an interesting case study for LLM reasoning due to its combination of mathematical theory and engineering design. We introduce ControlBench, a benchmark dataset tailored to reflect the breadth, depth, and complexity of classical control  design. We use this dataset to study and evaluate the problem-solving abilities of these LLMs in the context of control engineering. We present evaluations conducted by a panel of human experts, providing insights into the accuracy, reasoning, and explanatory prowess of LLMs in control engineering. Our  analysis reveals the strengths and limitations of each LLM in the context of classical control, and our results imply that Claude 3 Opus has become the state-of-the-art LLM for solving undergraduate control problems.  Our study serves as an initial step towards the broader goal of employing artificial general intelligence in control engineering. 
\end{abstract}

\section{Introduction}
Recently, the landscape of large language models (LLMs) has witnessed rapid advancements, with models such as GPT-4 \cite{achiam2023gpt}, Claude 3 \cite{Claude_3}, and Gemini 1.0 Ultra \cite{team2023gemini} pushing the boundaries of what artificial intelligence (AI) can achieve in complex problem-solving scenarios. 
These developments have sparked a growing interest in the application of LLMs across various domains, including coding \cite{nijkamp2022codegen, nam2024using, xu2022systematic, chew2023llm, macneil2022generating}, reasoning \cite{wei2022chain, huang2022towards,zhou2022least,sun2023survey, havrilla2024glore}, mathematics \cite{imani2023mathprompter, azerbayev2023llemma,frieder2024mathematical, zhang2024mathverse, he2023solving}, science \cite{wang2023scibench,birhane2023science,ouyang2023structured, yeadon2024impact, chen2023bioinfo}, and planning \cite{valmeekam2022large, valmeekam2024planning, zhao2024large, song2023llm, dagan2023dynamic}. 
Recent discussions and studies have highlighted the impressive capabilities of LLMs in many of these tasks. Following this trajectory, our paper aims to explore
the capabilities of state-of-the-art LLMs in solving undergraduate-level control system problems, a cornerstone of engineering education and research.

Automatic control is a fundamental pillar of modern engineering known for its complex feedback system designs and theoretical depth \cite{ogata2002modern,franklin1994feedback,nise2020control,aastrom2021feedback}.
The potential of LLMs to tackle control engineering problems presents an intriguing avenue for research, given the discipline's reliance on both mathematical rigor and engineering design. Control engineering encompasses a variety of challenging concepts including system dynamics, and PID / loopshaping design, and stability/robustness analysis of feedback mechanisms. The ability of LLMs to understand and solve undergraduate-level control problems could signify a substantial leap by integrating artificial intelligence into modern control engineering pipelines. More generally, the exploration of LLMs in solving control problems marks a significant inquiry into the potential of these models to contribute to areas traditionally reserved for specialized domain-specific human expertise. 

In this paper, we delve into this emerging research area by evaluating the capabilities of three state-of-the-art LLMs, namely GPT-4, Claude 3 Opus, and Gemini 1.0 Ultra, on our proposed college-level {\bf Control} system problem-solving {\bf Bench}mark, referred to as {\bf ControlBench}. Our proposed ControlBench
introduces a carefully crafted dataset designed to reflect the breadth, depth, and complexity of classical control design. Our problem set captures the essence of feedback system design and the analytical skills required in this field. Through a comprehensive evaluation conducted by a panel of human experts, we assess the performance of leading LLMs, including GPT-4, Claude 3 Opus, and Gemini 1.0 Ultra, in terms of their accuracy, reasoning capabilities, and ability to provide coherent and informative explanations.
Our analysis sheds light on the distinct strengths and limitations of each model, offering valuable insights into the potential role of LLMs in control engineering. This investigation not only contributes to our understanding of the current capabilities of LLMs but also paves the way for future research aimed at harnessing artificial general intelligence  in the advancement of control engineering solutions. 
Our main contributions can be summarized as follows.

\begin{itemize}
    \item  We introduce
a new natural-language dataset, called ControlBench, to test the capabilities
of LLMs in solving undergraduate control system problems.
    \item We present evaluations of GPT-4, Claude 3 Opus, and Gemini 1.0 Ultra on ControlBench, conducted by a panel of human experts. Built upon our accuracy and failure mode analysis, we further discuss the strengths and limitations of these LLMs. We present  various examples of LLM-based ControlBench responses  to support our discussion. Our results imply that Claude 3 Opus has become the state-of-the-art LLM in solving undergraduate control problems, outperforming the others in this study.  Based on our observation, one main limitation for all three LLMs is that they also suffer on problems involving visual elements such as Bode plots and Nyquist plots. Our study also sheds light on the role of self-correction, and issues such as sensitivity to the problem statements. 
    \item We also introduce a simplified version of ControlBench, termed as ControlBench-C, which only consists of single-answer multiple-choice problems. ControlBench-C enables fast automatic evaluation of LLMs from researchers without control background. We also highlight the limitations of ControlBench-C. Specifically, ControlBench-C is much simpler than ControlBench, and can not provide a comprehensive evaluation for the reasoning capabilities of LLMs in control engineering.
\end{itemize}
By examining the performance of LLMs in this specialized domain, we take an important step towards realizing the broader goal of integrating machine intelligence into the fabric of engineering education and research. 

{\bf Related Work:} There have been some early efforts on using LLMs to generate codes, cost functions, and/or value maps for robotic control tasks \cite{moweroptimal,huang2023voxposer}. Our paper takes a complementary angle, focusing on benchmarking the capabilities of various LLMs in solving undergraduate control problems.

\section{Motivating Example: A Showcase for LLM Capabilities in Control Design}
\label{sec:example1}

Before proceeding to a more comprehensive benchmark study, we will start with an illustrative example to showcase
 to showcase the potential of 
LLMs for solving control design problems. 
In this section,  we focus on a simple motivating example that considers  the design of a proportional-integral (PI) controller for a cruise control system.  We will ask GPT-4, Claude 3 Opus, and Gemini 1.0 Ultra to solve this simple PI design problem, and make a few observations. The exact problem statement is given as follows.\footnote{This example is actually Problem 5.17 from ControlBench, which will be introducec in the next section.}

\begin{conversationbox}
\centering{{\bf PI Design Example}}
    \begin{humaninput}
   \textbf{Problem Statement}: Consider a car whose longitudinal motion is modeled by the following ODE:
$$2085 \dot{v}(t) + 23.2v(t) = 40u(t) + 108.4 - F_{grav}(t)$$
The input is the throttle u and the output is the velocity v. The gravitational force $F_{grav}$ is a disturbance. Let $e(t) = v_{des}- v(t)$ denote the tracking error between the desired velocity $v_{des} = 29$ m/s and actual velocity v(t). Consider a PI controller of the following form:
$$ u(t) =  \bar{u} + K_p e(t) + K_i \int_0^t  e(\tau) d\tau  $$
where $\bar{u}=14.11$ is the open-loop input to maintain $v_{des}$ when on 
at road $\theta= 0^\circ$ . Choose the PI gains so that the cruise control system is stable and rejects disturbances due to changing road slopes within $\approx 10$ sec. The closed-loop should also be over or critically damped as oscillations are uncomfortable for the driver.
    \end{humaninput}
\end{conversationbox}

This particular problem is of interest as it is one of the first design problems encountered by undergraduates in  control engineering. Moreover, it requires reasoning about the impact of the control gains on the closed-loop poles and the connection to transient response properties. There is more than one design that satisfies the given objectives.  One approach is to note that the closed-loop characteristic equation is
$ 2085 s^2 + (23.2 + 40K_p) s + 40K_i =0$. 
We can place the closed-loop poles in the left half plane (LHP) to have any desired damping ratio $\zeta$ and natural frequency $\omega_n$ by choosing the controller gains to satisfy $\omega_n^2 = \frac{40K_i}{2085}$
and $2\zeta\omega_n = \frac{23.2+40K_p}{2085}$.  Selecting $\zeta=1$ gives two critically damped closed-loop poles placed at $s=-\omega_n$. The time constant for these poles is $\tau=\frac{1}{\omega_n}$ sec and, for critically damped poles, the 5\% settling time is approximately $4.75\tau = \frac{4.75}{\omega_n}$.  Thus we can select $\omega_n = \frac{4.75}{10} = 0.475$ to obtain a settling time near 10 sec.  Finally, using $\zeta=1$ and $\omega_n=0.475$ to solve for the PI gains leads to $K_p=48.94$ and $K_i=11.76$.  The corresponding transfer function from $F_{grav}$ to $v$ is
\begin{align}\label{eq:Tf}
    T_{Fgrav \to v}(s) = \frac{-s}{2085 s^2 + 1981 s + 470.4}
\end{align}
This has critically damped poles at $s=-0.475$ as expected. A step change increase in gravitational force (due to a step increase in the road slope) will cause the velocity to initially drop from the desired value.  However, the PI controller rejects the disturbance in $\approx 10$ sec with the velocity converging to the desired value with a nice overdamped response.

{\bf Initial Answers from LLMs:}
We tested GPT-4, Claude 3 Opus, and Gemini 1.0 on  the above problem. GPT-4 gives a short answer describing the process with a self-evaluated 40\% confidence in its answer (we ask LLMs to provide confidence scores for their own answers).  The key part of the GPT-4 answer makes the following statement without giving feasible values for the PI gains.

\begin{gptoutput}
\textbf{GPT-4}: Steps for PI Controller Design include Linearization, Transfer Function derivation, and PI Controller Design to ensure Stability, Disturbance Rejection, and Damping criteria are met.
\end{gptoutput}

\noindent We marked the answer from GPT-4 as "Incorrect". Gemini 1.0 Ultra gave a slightly more detailed set of design steps but also failed to provide values for the PI gains.  The solution from Gemini 1.0 Ultra is also marked as "Incorrect." 

Interestingly, in our initial testing, Claude 3 Opus followed almost the exact procedure as our baseline solution given above with a self-evaluated 90\% confidence.  It used a slightly different approximation for the settling time and chose closed-loop poles at $s=-0.1$ and $-0.2$ rad/sec.  It then solved for the gains $K_p=48.1$ and $K_i=10.4$.  It made a calculation error as these gains actually place the closed-loop poles at $s=-0.60$ and $-0.33$.  However, this calculation error is minor and the reasoning used by Claude 3 Opus in the design is correct.  One other minor point is that Claude 3 bases its design on the closed-loop transfer function from reference to output, i.e. on the complementary sensitivity.
\begin{align}
T(s) = \frac{40K_p s + 40K_i}{2085s^2 + (23.2 + 40K_p)s + 40K_i}
\end{align}
Note that this transfer function is distinct from the disturbance to output transfer function \eqref{eq:Tf} as it includes a zero  due to the PI controller. The gains given by Claude~3 Opus place this zero at at $-\frac{K_i}{K_p} = -0.21$. This zero is close to the imaginary axis and, as a result, causes a significant overshoot due to step changes in the reference speed.  This might be uncomfortable to the driver.  To investigate this issue further, we reworded the design problem to include the following prompt. 

\vspace{0.05in}
\emph{{\bf Reworded Problem:} To design a PI controller for the given car longitudinal motion model, we need to choose the gains
$K_p$ and $K_i$ such that the closed-loop system is stable, rejects disturbances within approximately 10
seconds, and has an overshoot less than 1\% in the step response.}
\vspace{0.05in}

This variation of the question is more challenging than the original wording as it requires the LLM to reason about the effect of the closed-loop zero on the step response characteristics. In other words, it is not sufficient to simply place overdamped closed-loop poles as the zero can causes overshoot as in the design solutions above.  Claude 3 Opus repeats its design with a confidence of 90\% and the following reasoning:

\begin{claude3output}
\textbf{Claude 3 Opus}: For a stable system with an overshoot less than 1\%, the damping ratio should be $\zeta \geq 0.826$. We can choose the poles to be $s_{1,2} = -\zeta\omega_n \pm j\omega_n\sqrt{1-\zeta^2}$, with $\zeta = 0.9$ and $\omega_n = 0.4$ to achieve a settling time of approximately 10 seconds (settling time $\approx \frac{4}{\zeta\omega_n}$).
\end{claude3output}

\noindent Claude 3 Opus then computes gains of $K_p=71.6$ and $K_i=33.4$. These gains yield closed-loop poles with $\zeta=0.865$ and $\omega_n=0.8$rad/sec so there again appears to be a minor calculation error.  More importantly, Claude 3 Opus failed to account for the effect of the zero on the overshoot in its reasoning.  Finally, we give the prompt: {\emph ``think about the impacts of closed-loop zeros."}
This prompt was sufficient as Claude 3 Opus then re-worked its design to account for the zero and returned PI gains $K_p=200$ and $K_i=5$ that did, indeed, give a step reference response with overshoot less than 1\%.  

Obviously, there is some randomness in the answer generation from LLMs. For the above question, even with the last prompt on zeros,  LLMs (including Claude 3 Opus) do not always give correct answers due to potential calculation mistakes. However, we believe that the above example demonstrates the potential of LLMs in control engineering, as they already have mastered the knowledge of basic control design to some extent. Next, we will introduce ControlBench for examining the capabilities of LLMs in control engineering.

\section{The ControlBench Dataset}
To assess the capabilities of LLMs in solving control problems,  we first create a collection of 147 undergraduate control problems. This collection comprises problems sourced from exercises in \cite{distefano1997schaum}, and problems gathered from undergraduate control classes at University of Michigan (EECS 460) and University of Illinois, Urbana-Champaign (ECE 486).

Our problem set spans a broad spectrum of topics typically encountered in undergraduate control courses, blending both textual and visual elements to mirror the multifaceted nature of real-world applications. Such integration is crucial as control system design inherently necessitates various types of plots to analyze and understand system behaviors. For instance, in the context of frequency-domain control design, Bode plots and Nyquist plots are often used as fundamental tools for analysis. Our dataset covers these topics, serving as a valuable tool for assessing the efficacy of LLMs in utilizing graphical information to tackle control problems. We summarize the statistics of our control problem dataset for each sub-topic in Table~\ref{tab:controlProblems}, where we also report the number of problems with visual elements under each topic.

We collect each problem from original documents in PDF files and presentation slides. We manually transfer these problems into LaTeX format. All the problems are carefully verified by human annotators to ensure that LaTeX
documents can be compiled without any syntax errors. In addition, we also provide a detailed step-by-step solution for each problem in LaTeX. Our ControlBench dataset and the related Latex/PDF files are available  at \url{https://agi4engineering.github.io/LLM4Control/}.

\begin{table}[t!]
\centering
\caption{Summary of the control problem dataset. We report the number of problems under each topic, and the number of problems with visual elements.}
\label{tab:controlProblems}
\begin{tabular}{lll}
\toprule
\textbf{Topic} & \textbf{\# of Problems}   & \textbf{\# Visual} \\ 
\midrule
Background & 28  &0   \\
Stability & 19  &0   \\
Time response & 21  &3   \\
Block diagrams & 5  &5   \\
Control System Design &24 &0 \\
Bode Analysis &15 &13 \\
Root-Locus Design &7 &1 \\
Nyquist Design &5 &4 \\
Gain/Phase Margins &9 &0 \\
System Sensitivity Measures  &3 &0 \\
Loop-shaping &4 & 0\\
Advanced Topics &7 &0 \\
\midrule 
Total  &147 &26  \\
\bottomrule
\end{tabular}
\end{table}

\section{Evaluations of Leading LLMs on ControlBench}\label{sec:eval}
In this section, we present the evaluation results for GPT-4, Claude 3 Opus, and Gemini 1.0 Ultra, focusing on their performance in solving control problems complied within our dataset ControlBench.

\subsection{Statistical Accuracy Analysis}

First, we examine the accuracy of GPT-4, Claude 3 Opus, and Gemini 1.0 Ultra in solving the problems from ControlBench.
To facilitate a foundational understanding of how LLMs approach control problems, we start with a zero-shot prompting strategy (i.e., directly inputting the question) and examine the responses through human annotation.  This setting can illustrated as follows:

\begin{conversationbox}
\centering{{\bf Zero-Shot Setting}}
    \begin{humaninput}
    \textbf{Human Input}: Consider the control problem \ldots
    \end{humaninput}
    \begin{gptoutput}
    \textbf{GPT-4}: For this problem, we need to consider several aspects, including the dynamics of the plant, the characteristics of the PI controller, and  \ldots
    \end{gptoutput}
     \begin{claude3output}
   \textbf{Claude 3 Opus}: For this problem, we need to consider the closed-loop system dynamics. Let's approach this step by step \ldots
    \end{claude3output}
    \begin{geminioutput}
      \textbf{Gemini 1.0 Ultra}:  Absolutely, let's analyze \ldots
    \end{geminioutput}
\end{conversationbox}

Once we get the LLM responses, we check the correctness of the LLM answers via human annotation. 
The input to LLMs is typically just copied from the LaTex description of the problems in ControlBench. 
This simple zero-shot setting serves a starting point for our analysis.
To understand how self-checking can improve the performance \cite{huang2022large,kojima2022large,wang2022self}, we also study a setting where a subsequent self-checking prompt such as "\textit{carefully check your answer again}" is used to aid LLMs in identifying and rectifying their previous mistakes. This setting is illustrated as follows.

\begin{conversationbox}
\centering{{\bf Self-Checking (or Self-Correction)}}
    \begin{humaninput}
    \textbf{Human Input }: Consider the control problem \ldots
    \end{humaninput}
    \begin{claude3output}
   \textbf{Claude 3 Opus}: We can perform the following calculations step by step \ldots
    \end{claude3output}

    \begin{humaninput}
    \textbf{Human Input}: Check your solutions carefully and fix all the errors.
    \end{humaninput}
    \begin{claude3output}
    \textbf{Claude 3 Opus}: I apologize for the confusion. Let me correct the errors and provide a more accurate solution \ldots
    \end{claude3output}
\end{conversationbox}

Once we get the revised LLM answers, the correctness is again examined by human experts. In this section, we test both settings for GPT-4, Claude 3 Opus, and Gemini~1.0~Ultra.

\textbf{Evaluation Metric:} Our main evaluation metric is Accuracy (\textbf{ACC}), defined as the proportion of instances where the LLMs correctly solve the given problems. In the LLM literature, it is known that sometimes a simple self-checking prompt, such as "\textit{carefully check your answer again}", enables LLMs to identity and fix errors in their previous answers \cite{huang2022large}. To quantify this self-correction ability, we deploy a subsequent self-checking prompt for initially incorrect responses, aiding LLMs in identifying and amending mistakes. This leads to the Self-Checked Accuracy (\textbf{ACC-s}) as a secondary metric, which quantifies the instances in which LLMs successfully amend their answers after a self-review process.

\begin{table*}[t]
  \centering
  \caption{Accuracy (ACC) and Self-Checked Accuracy (ACC-s) of GPT-4, Claude 3 Opus, and Gemini 1.0 Ultra in solving control problems across various topics. The ACC-s metric represents the models' ability to correct their initial errors upon a self-review, offering insight into the adaptability and error-correction capabilities of each model when tackling control problems. The best
results for each metric are highlighted in bold.}
  \label{tab: main result}
  \resizebox{1\columnwidth}{!}{%
  \begin{tabular}{l|cc|cc|cc}
    \toprule
     & \multicolumn{2}{c|}{\textbf{GPT-4}} & \multicolumn{2}{c|}{\textbf{Claude 3 Opus}} & \multicolumn{2}{c}{\textbf{Gemini 1.0 Ultra}} \\
   \midrule
   \textbf{Topics} & ACC $\uparrow$ & ACC-s $\uparrow$ & ACC & ACC-s & ACC & ACC-s \\
\midrule
Background & 60.7\%  (17/28)  &64.3\% (18/28)  & \textbf{75\% (21/28)} &\textbf{89.3\% (25/28)}  & 53.6\% (15/28)  &57.1\% (16/28) \\
Stability & 57.9\% (11/19)  &57.9\% (11/19)  & \textbf{76.2\% (15/19)}  &\textbf{89.5\% (17/19)}  & 31.6\% (6/19)  &31.6\% (6/19)  \\
Time response & 57.1\% (12/21) &66.6\% (14/21) & \textbf{76.2\% (16/21)}  &\textbf{76.2\% (16/21)}  & 52.4\% (11/21) &57.1\% (12/21)  \\
Block diagrams &\textbf{40.0\% (2/5)} &\textbf{40.0\% (2/5)} &\textbf{40.0\% (2/5)}  &\textbf{60.0\% (3/5)}  & 0.0\% (0/5)  &0.0\% (0/5)  \\
Control System Design &29.2\% (7/24) &29.2\% (7/24) &\textbf{ 33.3\% (8/24)}  &\textbf{62.5\% (15/24)}  & 25.0\% (6/24)  &37.5\% (9/24) \\
Bode Analysis &6.66\% (1/15) &6.66\% (1/15) & \textbf{13.3\% (2/15)}  &\textbf{13.3\% (2/15)}  & 6.66\% (1/15)  &6.66\% (1/15)\\
Root-Locus Design &28.6\% (2/7) &28.6\% (2/7) & \textbf{42.9\% (3/7)} &\textbf{42.9\% (3/7)}  & 28.6\% (2/7)  &28.6\% (2/7) \\
Nyquist Design &0.0\% (0/5) &0.0\% (0/5) &\textbf{40.0\% (2/5)}  &\textbf{40.0\% (2/5)} & 0.0\% (0/5)  &0.0\%  (0/5)\\
Gain/Phase Margins &\textbf{66.7\% (6/9)} &\textbf{66.7\% (6/9)} & \textbf{66.7\% (6/9)}  &\textbf{66.7\% (6/9)}  & 33.3\% (3/9)  &33.3\% (3/9) \\
System Sensitivity Measures  &\textbf{100.0\% (3/3)} &\textbf{100.0\% (3/3)} & \textbf{100.0\% (3/3)}  &\textbf{100.0\% (3/3)}  & 66.7\% (2/3)  &\textbf{100.0\% (3/3)} \\
Loop-shaping  &25.0\% (1/4) &25.0\% (1/4) & \textbf{50.0\% (2/4)}  &\textbf{75.0\% (3/4)}  & 25.0\% (1/4)  &25.0\% (1/4) \\
Advanced Topics &71.4\% (5/7) &71.4\% (5/7) & \textbf{85.7\% (6/7)}  &\textbf{85.7\% (6/7)}  & 42.9\% (3/7)  &57.1\% (4/7) \\
\midrule 
Total &45.6\% (67/147) &47.6\% (70/147) & \textbf{58.5\% (86/147)}  &\textbf{68.7\% (101/147)}  & 34.0\% (50/147)  &38.8\% (57/147)  \\ 
\bottomrule
\end{tabular}
}
\end{table*}
Table \ref{tab: main result} provides the Accuracy (ACC) and Self-Checked Accuracy (ACC-s) of GPT-4, Claude 3 Opus, and Gemini 1.0 Ultra across a variety of control problem topics in our dataset ControlBench. Claude 3 Opus emerges as the standout model, demonstrating superior performance in both ACC and ACC-s. This indicates that Claude 3 Opus not only has a higher baseline understanding of control problems but also exhibits superior self-correction capabilities\footnote{This aligns with the findings in \cite{Claude_3} mentioning the superior performance of Claude 3 Opus on many tasks compared to other LLMs including GPT-4 and Gemini 1.0 Ultra.}. 
GPT-4 displays competitive accuracy in specific areas, such as Block Diagrams, Root-Locus Design, and System Sensitivity Measures, but does not match Claude 3 Opus in overall performance. Gemini 1.0 Ultra, meanwhile, ranks lower in general effectiveness. Notably, Claude 3 Opus also shows enhanced proficiency in handling problems with visual components, particularly in topics like Nyquist Design, Root-Locus Design, and Loop-shaping.
It is important to clarify that the accuracies detailed in Table \ref{tab: main result} are not intended for direct comparison with human expertise. Instead, these results establish an initial benchmark for the performance of LLMs in addressing college-level control problems.

We also want to emphasize that there is some randomness in the answer generations of all three LLMs. Even for the same question, all three LLMs can generate drastically different answers. Therefore, our results should be interpreted in a more qualitative manner, i.e. they give an initial comparison of the capabilities 
of GPT-4, Claude 3 Opus, and Gemini 1.0 Ultra in solving control problems. More study is definitely needed in the future. Next, we showcase the potential of LLMs in solving control design problems via an illustrative example, and then discuss their failure modes and self-correction capabilities observed in our experiments.

\subsection{Strengths of LLMs and Successful Examples}
\label{sec:success}

From Table \ref{tab: main result}, we can see that LLMs (especially Claude 3 Opus and GPT-4) do reasonably well in solving problems related to background mathematics, stability, and time response. In addition, Claude 3 Opus has achieved $33.3\%$ ACC and $70.8\%$ ACC-s in basic control design problems, demonstrating its potential in the domain of control system design. In this section, we will present several successful examples to showcase the current capabilities of LLMs in control engineering.  We will start with the following simple stability problem which has been succesfully solved by all three LLMs.

\begin{conversationbox}
\centering{{\bf Routh-Hurwitz Stability Criteria (Problem 2.6 from ControlBench)}}
    \begin{humaninput}
    \textbf{Human Input}: A system has a characteristic equation given by
\begin{equation}
    s^3 + Ks^2 + (K+1)s + 6 = 0.
\end{equation}
Determine the range of \(K\) for a stable system. You can use Routh-Hurwitz Stability Criteria.
    \end{humaninput}
    \begin{gptoutput}
    \textbf{GPT-4}: To determine the range of \(K\) for a stable system using the Routh-Hurwitz Stability Criteria, we start with the given characteristic equation:
\[
s^3 + Ks^2 + (K+1)s + 6 = 0.
\]

We construct the Routh array as follows \ldots 
Solving the second inequality, we find the solution to be \(K < -3\) or \(K > 2\). However, since \(K > 0\) is required for stability, we ignore the solution \(K < -3\).
Hence, for the system to be stable, the range of \(K\) is
$K > 2$.
\end{gptoutput}

\begin{claude3output}
\textbf{Claude 3 Opus}: To determine the range of K for a stable system using the Routh-Hurwitz stability criteria, we need to analyze the coefficients of the characteristic equation.
Given:
The characteristic equation is $s^3 + Ks^2 + (K+1)s + 6 = 0$.
Step 1: Write the characteristic equation in the standard form.
$a_0s^3 + a_1s^2 + a_2s + a_3 = 0$, where:
$a_0 = 1$,
$a_1 = K$,
$a_2 = K+1$,
$a_3 = 6$.
Step 2: Apply the Routh-Hurwitz stability criteria for a third-order system.
For a stable system, all the coefficients of the characteristic equation must be positive, and the following condition must be satisfied:
$a_1a_2 > a_0a_3$.
Step 3: Check the positivity of the coefficients.
$a_0 = 1 > 0$ (always positive)
$a_1 = K > 0$
$a_2 = K+1 > 0$
$a_3 = 6 > 0$ (always positive)
From $a_1 > 0$ and $a_2 > 0$, we get:
$K > 0$ and $K > -1$.
Step 4: Check the condition $a_1a_2 > a_0a_3$.
$a_1a_2 > a_0a_3$
$K(K+1) > 1 \times 6$
$K^2 + K - 6 > 0$.
Solving the quadratic inequality:
$(K + 3)(K - 2) > 0$
$K < -3$ or $K > 2$.
Step 5: Combine the conditions for K.
From steps 3 and 4, we have:
$K > 0$ and $K > 2$.
Therefore, the range of K for a stable system is $K > 2$.
\end{claude3output}
\begin{geminioutput}
    \textbf{Gemini 1.0 Ultra}: Routh-Hurwitz Stability Criterion: The Routh array is:
    
$\begin{array}{c|cc}
s^3 & 1 & K+1 \\
s^2 & K & 6 \\
s^1 & \frac{K^2 + K - 6}{K} & 0 \\
s^0 & 6 & 0
\end{array}$

Conditions for Stability: $K > 0$, $K^2 + K - 6 > 0 \implies (K-2)(K+3) > 0$, $6 > 0$ (always true)

Range of K for stability: $K > 2$
\end{geminioutput}
\end{conversationbox}

\noindent For the above stability problem, all three LLMs follow roughly the same reasoning path and reach to the correct answer eventually. There are many other stability problems (from ControlBench) which actually cause difficulty for Gemini 1.0 Ultra or even GPT-4.  Overall, we can see from
 Table \ref{tab: main result} that Claude 3 Opus  has outperformed GPT-4 and Gemini 1.0 Ultra in answering ControlBench questions related to background mathematics, stability, and time response.

Next, we give one example showcasing the basic control design knowledge of GPT-4 and Claude 3 Opus. 

\begin{conversationbox}
\centering{{\bf Pole Placement (Problem 5.2 from ControlBench)}}
    \begin{humaninput}
    \textbf{Human Input}: A plant with nominal model 
\begin{equation}
    G_o(s) = \frac{1}{(s + 1)^2}
\end{equation}
is in a feedback loop under control with a PI controller having transfer function $C(s) = \frac{K_ps + K_i}{s}$.
Determine whether this controller configuration can be used to achieve full pole placement.
    \end{humaninput}
    \begin{gptoutput}
    \textbf{GPT-4}: In control theory, pole placement is a strategy for the design of feedback controllers to force the closed-loop control system to have desired eigenvalues or poles \ldots 
    This characteristic equation is a cubic polynomial in 
$s$, which means it will have three roots or poles. However, the ability to freely place these poles is limited by the structure of the polynomial. In general, full pole placement requires as many adjustable parameters as there are poles. In this case, we have only two parameters ($K_p$ and $K_i$) to adjust three poles.
Thus, while the PI controller configuration can influence the location of the poles to some extent, it cannot achieve arbitrary pole placement for all three poles due to the limited number of tuning parameters. Additional control terms, such as a derivative term (making it a PID controller), may be required to achieve full pole placement capability in this scenario.
\end{gptoutput}

\begin{claude3output}
\textbf{Claude 3 Opus}: To determine if the given PI controller configuration can achieve full pole placement, we need to analyze the closed-loop transfer function of the system \ldots 
For full pole placement, the controller should have at least as many parameters as the order of the plant. In this case, the plant is of second order, and the PI controller has two parameters ($K_p$ and $K_i$). However, the closed-loop characteristic equation is of third order.
Therefore, the given PI controller configuration cannot achieve full pole placement for this plant. To achieve full pole placement, we would need a controller with at least three parameters, such as a PID controller or a higher-order controller.
\end{claude3output}
\begin{geminioutput}
    \textbf{Gemini 1.0 Ultra}: Here's a breakdown of how to determine if full pole placement is achievable  \ldots Because the number of adjustable parameters from the controller matches the number of poles in the plant, you can likely achieve full pole placement with the given PI controller configuration. Note that this assumes the system is controllable; you'd still need to verify its controllability.
\end{geminioutput}
\end{conversationbox}

\noindent
For the above problem, both GPT-4 and Claude 3 Opus reach 
the correct conclusion, and the sampled answer from Gemini 1.0 Ultra is wrong. If we keep on sampling answers for this question, sometimes Gemini 1.0 Ultra will also give the correct answer (and Claude 3 Opus can occasionally give a wrong answer at first and then fixes the mistake after self-checking).  This example gives us a rough idea on the knowledge level of LLMs on basic control design. 
Based on Table \ref{tab: main result}, it is fair to state that Claude 3 Opus has shown stronger potential in addressing classic control design tasks, i.e. it has achieved much higher ACC and ACC-s for this topic. Now we give another example to showcase the superior design capabilities of Claude 3 Opus.

\begin{conversationbox}
\centering{{\bf PID Control (Problem 5.11 from ControlBench)}}
    \begin{humaninput}
    \textbf{Human Input}: Consider the plant with the following transfer function:
\begin{align*}
    G(s) = \frac{505}{s^3+21s^2+121s+101}
\end{align*}
The poles of this system are at \(s = -1, -10 \pm j\), and the DC gain is $5$.
\begin{itemize}
    \item[(a)] What is the dominant pole approximation \(G_a(s)\) for this plant?
    \item[(b)] Would you recommend using a PI, PD, or PID controller?
    \item[(c)] Choose the controller gains so that the closed-loop with \(G_a(s)\) has poles repeated at \(s = -1\)
\end{itemize}
    \end{humaninput}
    \begin{gptoutput}
    \textbf{GPT-4}: To address your questions, let's break down each part systematically \ldots Therefore, the dominant pole approximation $G_a(s)$ is $G_a(s)=\frac{5}{s+1}$ \ldots  
Given the system has a single dominant pole and assuming we want both good steady-state performance and responsive dynamics, a PID controller would be a balanced choice, providing a compromise between responsiveness and steady-state error correction \ldots 
Given the complexity of 
$G(s)$ and the nonlinearity introduced by the PID controller, this process can be quite involved and may require numerical methods or root locus techniques to solve accurately. If a more detailed solution is needed, I could assist further with that process.
\end{gptoutput}

\begin{claude3output}
\textbf{Claude 3 Opus}: Alright, let's tackle this problem step by step:

(a) The dominant pole approximation $G_a(s)$ is found by considering the pole closest to the origin (the dominant pole) and neglecting the other poles. In this case, the dominant pole is at $s = -1$. The transfer function of the dominant pole approximation is:
$G_a(s) = \frac{K}{s+1}$,
where $K$ is the DC gain of the original system. Given that the DC gain is 5, we have:
$G_a(s) = \frac{5}{s+1}$.

(b) Since the dominant pole is real and the system is stable, a PI controller would be a good choice to eliminate steady-state error and improve the system's response.

(c) To place the closed-loop poles of $G_a(s)$ with a PI controller at $s = -1$, the characteristic equation should be:
$(s+1)^2 = 0$.
The closed-loop transfer function with a PI controller is:
$\frac{G_a(s)C(s)}{1+G_a(s)C(s)}$, where $C(s) = K_p + \frac{K_i}{s}$.
The characteristic equation is:
$1+G_a(s)C(s) = 0$,
$1+\frac{5}{s+1}(K_p + \frac{K_i}{s}) = 0$,
$1+\frac{5K_ps+5K_i}{s(s+1)} = 0$,
$s(s+1)+5K_ps+5K_i = 0$,
$s^2+(5K_p+1)s+5K_i = 0$.
To place the poles at $s = -1$, the characteristic equation should be:
$(s+1)^2 = s^2+2s+1 = 0$.
Comparing the coefficients, we get:
$5K_p+1 = 2$ and $5K_i = 1$.
Solving for the controller gains:
$K_p = \frac{1}{5}$ and $K_i = \frac{1}{5}$
Therefore, the PI controller transfer function is:
$C(s) = \frac{1}{5} + \frac{1}{5s}$.
\end{claude3output}
\begin{geminioutput}
    \textbf{Gemini 1.0 Ultra}: Absolutely, let's analyze this control system problem \ldots 
    We can approximate the transfer function as:
$G_a(s) \approx  \frac{505}{s^2 + 20s + 100}$\ldots
    We need to substitute the plant and controller transfer functions into the characteristic equation and equate coefficients. This is a somewhat involved process that would typically make use of symbolic computation software \ldots
\end{geminioutput}
\end{conversationbox}

\noindent For this problem, the above sampled answer from Claude 3 Opus is correct\footnote{Again, there is some randomness involved. Based on our observation on this problem, more often than not, Claude 3 gives the right answer.}. Interestingly, for Part (b), Claude 3 Opus often prefers a PI controller, while GPT-4 and Gemini 1.0 Ultra tend to choose  a PID controller (we have tried generating many sampled answers from GPT-4, and we observed that very occasionally, GPT-4 may also choose a PI controller).  Claude 3 Opus gives the following reason: `` a PI controller can effectively eliminate steady-state error and achieve the desired closed-loop performance without the added complexity of the derivative action." This captures the essence of this question: derivative control is not required as the dominant pole approximation for the plant is first order. It is interesting that Claude 3 Opus is able to correctly reason about the impact of this dominant pole approximation on the choice of the control architecture.
 This example again showcases the potential of Claude 3 Opus, although the consistency is still an open issue.

\subsection{Analysis and Insights for Failure Modes}
\label{sec:FM}

Despite their great potential, LLMs can fail in many different ways. 
In this section, we discuss various failure modes observed in our experiments, subsequently building towards more intricate failure mode analysis. 

\begin{figure}[t!]
    \centering
    \includegraphics[width=0.6\textwidth]{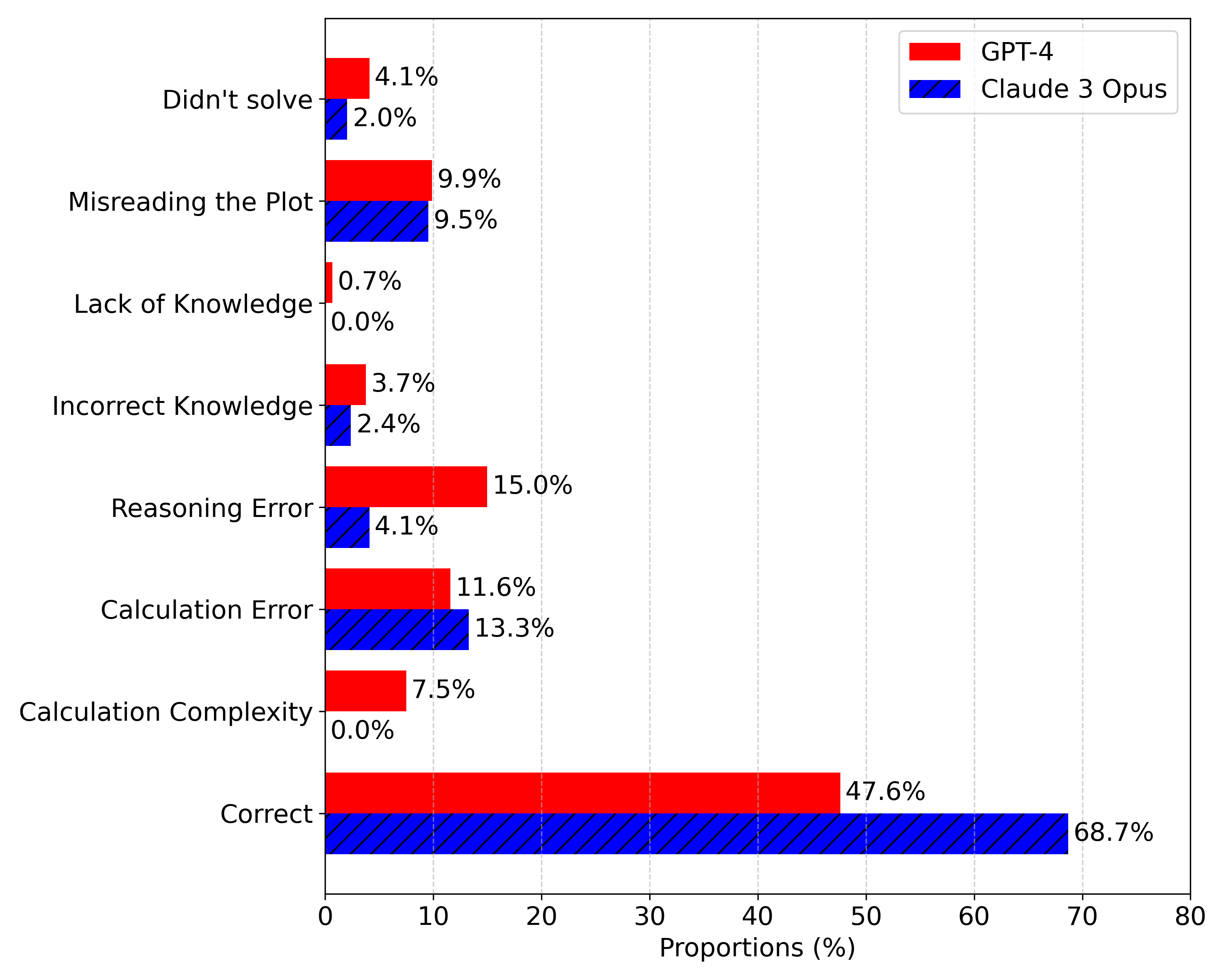}
    \caption{Proportions (\%) of seven error types (\#errors / \#total cases) for GPT-4 and Claude 3 Opus, including portion of correct answers for quick comparison}

     \label{fig:errortype}
    \vspace{-0.2in}
\end{figure}

We have observed quite diverse failure modes in our experiments. 
To make the discussion more formal, we categorize the LLM failures into  several types, and highlight the proportions (\%) of seven error types for GPT-4 and Claude 3 Opus in Figure \ref{fig:errortype} (for simplicity, our discussion in this section mainly focuses on comparing GPT-4 and Claude 3 Opus, both of which outperform Gemini 1.0 Ultra). Our error analysis is performed for ACC-s. From Figure~\ref{fig:errortype}, we can see that the biggest bottleneck preventing GPT-4 to achieve better accuracy on ControlBench is its limited reasoning capabilities. In contrast, Claude 3 Opus has significantly fewer  reasoning errors. This observation is consistent with the fact that Claude 3 Opus has surpassed GPT-4 for many reasoning tasks in general. Currently, the biggest bottleneck for Claude 3 Opus is its calculation abilities. However, issues related to calculation errors can be typically mitigated by calling external calculation tools or adopting the program-of-thought approach \cite{chen2022program}. For example, as demonstrated in \cite{ouyang2023structured}, the program-of-thought approach can significantly reduce the LLM calculation errors in chemistry problems via translating the LLM solution into \texttt{Python} codes. Although pursuing such extensions is beyond the scope of our paper, we expect that the performance of Claude 3 Opus on ControlBench can be significantly improved in the near future via improving the calculation precision.  In comparison to the calculation errors, the reasoning issue is more fundamental. In the machine learning field, there are on-going research efforts aiming to improve the reasoning capabilities of LLMs \cite{yao2024tree,besta2023graph,hao2023reasoning}. We will investigate how to integrate such efforts with ControlBench in the future.

We also notice that all three LLMs suffer from misreading graphical data including Bode, Nyquist, and root-locus plots. To illustrate this issue, we consider a simple Bode plot problem taken from EECS 460 at University of Michigan. The problem statement is given as follows (this is Problem 6.3 of ControlBench).

\vspace{0.1in}

 \begin{figure}[h]
\centering
\includegraphics[width=0.45\linewidth]{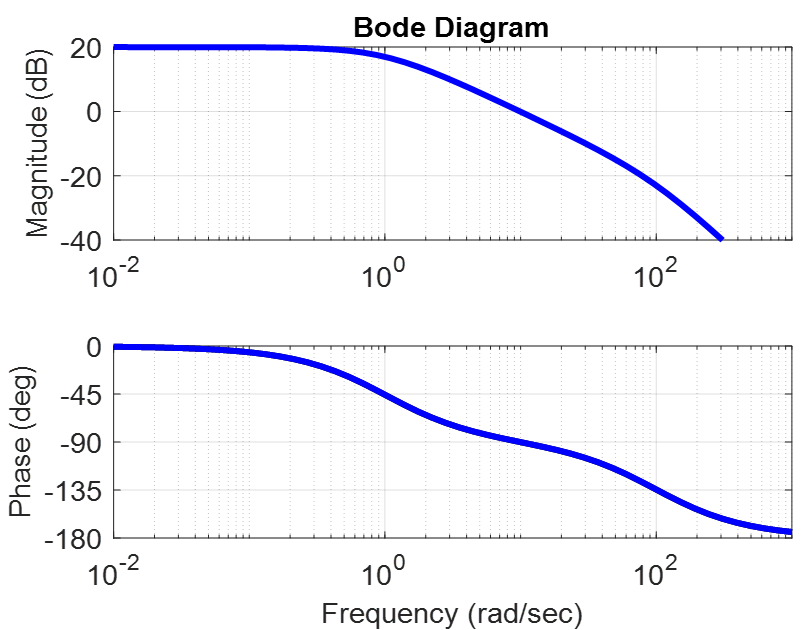}
\caption{Bode Plot Example}
\label{fig:bode}
\end{figure}

\begin{conversationbox}
\centering{{\bf Bode Plot Example (Problem 6.3 from ControlBench)}}
    \begin{humaninput}
    \textbf{Human Input}: A linear time-invariant system \(G(s)\) with input \(u\) and output \(y\) has the Bode plot shown above:
\begin{itemize}
    \item[(a)] What is \(|G(10j)|\) in \(dB\) and actual units?
    \item[(b)] What is \(\angle G(10j)\) in degs and radians?
\end{itemize}
    \end{humaninput}
    \begin{gptoutput}
    \textbf{GPT-4}: For (a), \( |G(10j)| \) in dB is approximately -6 dB. In actual units, \( |G(10j)| = 10^{\frac{-6}{20}} \approx 0.5012 \). For (b), \( \angle G(10j) \) is approximately -45 degrees, or in radians, \( \frac{-45 \times \pi}{180} \approx -0.7854 \) radians.
\end{gptoutput}

\begin{claude3output}
\textbf{Claude 3 Opus}: For (a), from the magnitude plot, $|G(10j)| \approx -45$ dB. Converting to actual units:
$|G(10j)| = 10^{-45/20} \approx 0.00562$. For (b), from the phase plot, $\angle G(10j) \approx -90^\circ$. Converting to radians: 
$\angle G(10j) \approx -\frac{\pi}{2}$ rad.
\end{claude3output}
\end{conversationbox}

\noindent
As humans, we can easily see that \(|G(j10)| = 0 dB = 1\) and \(\angle G(j10) = -90 \text{deg} = -\frac{\pi}{2} \text{rad}\). However, for such a simple task, GPT-4 and Claude 3 Opus both give wrong answers.
 Both LLMs misread the magnitude plot in Figure \ref{fig:bode}, and made mistakes in answering the first part of the above question. 
For this example, Claude 3 Opus somehow reads the phase plot correctly. However, in general, both GPT-4 and Claude 3 Opus perform badly on problems that require extracting data from plots. This highlights the need of developing more advanced visual language models for solving such problems. 
Besides the above main observations, several other findings are also worth discussing.  Next we present these findings.

\textbf{Common vs. Isolated Mistakes:} We have observed that GPT-4, Claude 3 Opus, and Gemini 1.0 Ultra sometimes make the same mistakes. However, it is more often to see that these models exhibit diverse isolated failure cases that are due to either the randomness in the LLM answer generation or
some specific training data idiosyncrasies. We illustrate this via the following simple example (Problem 3.5 from ControlBench).

\begin{conversationbox}
\centering{{\bf Problem 3.5 from ControlBench}}
    \begin{humaninput}
        \textbf{Human Input}: What is the unit step response of a continuous system whose transfer function has a zero at \(-1\), a pole at \(-2\), and a gain factor of \(2\)?
    \end{humaninput}
    \begin{gptoutput}
        \textbf{GPT-4}: Given a system with a transfer function characterized by a zero at \(s = -1\), a pole at \(s = -2\), and a gain factor of 2, the transfer function \(H(s)\) is expressed as $H(s) = 2\frac{s + 1}{s + 2}$.
The unit step response of the system can be found by determining the output \(Y(s)\) when the input is a unit step function \(u(t)\), whose Laplace transform is \(\frac{1}{s}\). Thus, the Laplace transform of the output is given by \(Y(s) = H(s) \times \frac{1}{s}= = 2\frac{s + 1}{s(s + 2)}\).
This expression can be decomposed into partial fractions to facilitate the inverse Laplace transform 
$Y(s) = 2\left( \frac{A}{s} + \frac{B}{s + 2} \right)$.
Solving for coefficients \(A\) and \(B\) by setting \(s=0\) and \(s=-2\), respectively, and then equating coefficients, we find
$Y(s) = \frac{1}{s + 2} + \frac{1}{s}$.
Taking the inverse Laplace transform of \(Y(s)\) gives the unit step response \(y(t)\) in the time domain:
$y(t) = \mathcal{L}^{-1}\left\{ \frac{1}{s + 2} + \frac{1}{s} \right\} = e^{-2t}u(t) + u(t)$,
where \(u(t)\) is the unit step function. Therefore, the unit step response of the given system is:
$y(t) = e^{-2t} + 1 \quad \text{for} \quad t \geq 0$.

    \end{gptoutput}
    \begin{claude3output}
    \textbf{Claude 3 Opus}: 
To find the unit step response of the given system, we need to follow these steps:

Write the transfer function of the system.
Perform partial fraction expansion on the transfer function.
Find the inverse Laplace transform of the resulting expression.

Step 1: Transfer function
The transfer function is given by:
$G(s) = 2 \cdot \frac{s+1}{s+2}$.

Step 2: Partial fraction expansion
$\frac{2(s+1)}{s+2} = 2 - \frac{2}{s+2}$

Step 3: Inverse Laplace transform
The inverse Laplace transform of a constant is the constant itself multiplied by the unit step function, and the inverse Laplace transform of \(\frac{1}{s+a}\) is \(e^{-at}u(t)\), where \(u(t)\) is the unit step function. Therefore, the inverse Laplace transform of \(2 - \frac{2}{s+2}\) is:
$$\mathcal{L}^{-1}\left[2 - \frac{2}{s+2}\right] = 2u(t) - 2e^{-2t}u(t)$$

The unit step response of the system is:
$y(t) = 2u(t) - 2e^{-2t}u(t)$.
This can be simplified to:
$y(t) = 2(1 - e^{-2t})u(t)$.
Therefore, the unit step response of the continuous system with a zero at \(-1\), a pole at \(-2\), and a gain factor of \(2\) is \(y(t) = 2(1 - e^{-2t})u(t)\).
    \end{claude3output}
\end{conversationbox}

\noindent For  the above problem, GPT-4 is correct, while Claude~3 Opus makes a mistake. In fact, no matter how we ask Claude 3 Opus to do self-checking on this problem, it almost always gives the same wrong answer here (interestingly, different wrong reasoning paths show up to lead to the same wrong final answer). 
Therefore, it seems that Claude 3 Opus has a special isolated failure mode for this problem, which is not shared by GPT-4 or Gemini 1.0 Ultra. 
Here is a sampled answer from Claude 3 Opus after self-checking.

\begin{conversationbox}
\centering{{\bf Claude 3 Opus with Self-Correction for Problem 3.5 from ControlBench}}
    \begin{humaninput}
        \textbf{Human Input}: Your answer is incorrect. Fix it.
    \end{humaninput}
    \begin{claude3output} \textbf{Claude 3 Opus}: I apologize for the confusion in my previous response. Let me correct the answer for you.
Given zero at \(s=-1\),
pole at \(s=-2\), gain factor of \(2\), the correct transfer function is:
$$G(s) = 2 \cdot \frac{s+1}{s+2}$$

To find the unit step response, we can use the following property of the Laplace transform:
$\mathcal{L}{u(t)} = \frac{1}{s}$.
Multiplying the transfer function by \(\frac{1}{s}\):
$$\frac{Y(s)}{U(s)} = G(s) = 2 \cdot \frac{s+1}{s+2} \cdot \frac{1}{s}$$

Simplifying:
$$Y(s) = 2 \cdot \frac{s+1}{s(s+2)}$$

Using partial fraction expansion:
$$Y(s) = 2 \cdot \left(\frac{1}{s} - \frac{1}{s+2}\right)$$

Taking the inverse Laplace transform:
$$y(t) = 2 \cdot (1 - e^{-2t})u(t)$$

Therefore, the correct unit step response of the continuous system with a zero at \(-1\), a pole at \(-2\), and a gain factor of \(2\) is \(y(t) = 2 \cdot (1 - e^{-2t})u(t)\).
  \end{claude3output}
\end{conversationbox}

\noindent
It is interesting to see how Claude 3 Opus gets stuck with the same mistake for this problem. 
If we prompt Claude 3 Opus to carefully check the partial fraction expansion of $Y(s)$, then it will generate the correct answer and give a very detailed reasoning process for getting the correct inverse Laplace transform of \(Y(s) = 2 \cdot \frac{s+1}{s(s+2)}\). This example demonstrates that each LLM may have its own internal bias which may affect the final solution and lead to unique isolated mistakes.

\textbf{Correlations of Failures and Confidence:} We can ask LLMs to provide a confidence score (\%) for their own answers. We observe that it is more likely that LLMs will give correct answers when their confidence level is high. 
However, we also observe that LLMs can give wrong answers even if the self-evaluated confidence level is high. For instance, for the above problem (3.5 from ControlBench), Claude 3 Opus is 95\% confident in its solution \(y(t) = 2 - 2e^{-2t}\), which turns out to be wrong. This is quite common for all three LLMs. Therefore, high confidence provided by LLMs does not mean high accuracy.

\textbf{Failures from Mathematical Derivations:} Another observation is that LLMs can give unexpected mathematical derivation errors,  which we believe that could be fixed using external symbolic tools such as \texttt{Mathematica}. We illustrate this on the following simple example problem.

\begin{conversationbox}
\centering{{\bf Problem 1.3 from ControlBench}}
    \begin{humaninput}
        \textbf{Human Input}: Consider a system with a transfer function 
$H(s) = \frac{2}{s + 1} + \frac{\alpha}{s + 2}$,
where \(\alpha\) is a real number. Is there a range of real values for \(\alpha\) such that the system's unit step response exhibits~undershoot?
    \end{humaninput}
    \end{conversationbox}

\noindent For this problem, it is sufficient  combine the terms of $H(s)$ and derive that it has a Non-Minimum Phase (NMP) zero, i.e. RHP zero, when  $-4<\alpha<-2$.
With prompts related to NMP zeros and self-checking, both GPT-4 and Claude 3 Opus  give the condition $-\frac{4+\alpha}{2+\alpha}>0$. However, neither is able to convert this condition to obtain the final correct answer $-4<\alpha<-2$. This example may reveal that LLMs have difficulties dealing with multiple symbolic inequalities. There are actually many more examples showing that LLMs can make unexpected errors in symbolic manipulations (see our project website for more examples).
This may motivate the future study on integrating LLMs with symbolic toolboxes for improving performances on ControlBench.

\subsection{Discussions on Self-Correction Capabilities}

In this section, we want to highlight the self-correction capabilities of LLMs \cite{huang2022large,gou2023critic}. It seems that Claude 3 Opus significantly improved the accuracy  after we prompt it to check its previous wrong solutions.    GPT-4 and Gemini 1.0 Ultra have also shown some improvements with the self-correction step, but the improvements are less significant.
For illustrations, consider the following simple stability problem.

\begin{conversationbox}
\centering{{\bf  Problem 2.1 from ControlBench}}
    \begin{humaninput}
    \textbf{Human Input}: In a feedback control loop, the open-loop transfer function \(G(s)\) and the controller \(C(s)\) are given by
\begin{align*}
G(s) = \frac{s - 2}{(s - 1)(s + 4)}, \quad C(s) = K \frac{s + 1}{s}
\end{align*}
Determine \(K \in \mathbb{R}\), if exists, such that the control loop is stable.
    \end{humaninput}
    \begin{claude3output}
    \textbf{Claude 3 Opus}: To determine the range of values for the gain $K$ that makes the closed-loop system stable, we need to analyze the characteristic equation of the system using the Routh-Hurwitz stability criterion \ldots Combining the conditions, we find that for the system to be stable, $K$ must satisfy:
$-3 < K < 2$.
Therefore, the range of values for $K$ that makes the closed-loop system stable is $(-3, 2)$.
    \end{claude3output}

    \begin{humaninput}
    \textbf{Human Input}: Carefully check your solution.
    \end{humaninput}
    \begin{claude3output}
    \textbf{Claude 3 Opus}: Thank you for the suggestion to double-check the solution. Upon reviewing the work, I found an error in the simplification of the characteristic equation. Let me correct the mistake and solve the problem again...We find that there is no value of $K$ that satisfies all three conditions simultaneously. Therefore, the closed-loop system cannot be made stable for any real value of $K$.
I apologize for the earlier mistake and thank you for prompting me to review the solution carefully.
    \end{claude3output}
\end{conversationbox}

The correct answer for the above problem is that there is no such real value for $K$ to make the system stable. At the beginning, Claude 3 Opus made a mistake in the calculations. Similarly, GPT-4 and Gemini 1.0 Ultra also do not give the correct answer. Surprisingly, after asking Claude 3 Opus to do self-checking, it does state that no such values of $K$ exists\footnote{We observed that due to the randomness in the answer generation, Claude 3 Opus sometimes corrects its mistake for this problem, and obtained the right final answer using only a partially correct reasoning path.}.
As a matter of fact, the self-correction ability of Claude 3 Opus is quite noticeable, i.e. it improves the accuracy on ControlBench by $13.6\%$  after self-correction.
We provide another interesting example to illustrate this.

\begin{conversationbox}
\centering{{\bf Claude 3 Opus with Self-Correction on PID Design (Problem 5.13 from ControlBench)}}
    \begin{humaninput}
    \textbf{Human Input}: Consider the plant with the following transfer function:
    $G(s) = \frac{20}{s^2-6s+10}$.
\begin{itemize}
    \item[(a)] What is the closed-loop ODE from reference \(r\) to output \(y\) if we use the following PID controller?
    $u(t) = K_p e(t) + K_i \int e(t) + K_d \dot{e}(t)$,
    where \(e(t) = r(t) - y(t)\)
    \item[(b)] Choose the controller gains so that the closed-loop has poles repeated at \(s = -3\). Hint: \((s+3)^3 = s^3+9s^2+27s+27\)
    \item[(c)] What is the impact of the implementing the derivative term \(K_d \dot{e}(t)\) as versus the rate feedback from \(-K_d \dot{y}(t)\)?
\end{itemize}

    \end{humaninput}
    \begin{claude3output}
    \textbf{Claude 3 Opus}: Let's tackle this problem step by step:

(a) The closed-loop ODE from reference $r$ to output $y$ using the PID controller can be derived as follows \ldots
Taking the inverse Laplace transform, we get the closed-loop ODE:

   $\dddot{y} + (20K_d-6)\ddot{y} + (20K_p+10)\dot{y} + 20K_i y = 20(K_d \ddot{r} + K_p \dot{r} + K_i r)$

(b) To choose the controller gains so that the closed-loop has poles repeated at $s = -3$, we compare the denominator of the closed-loop transfer function with $(s+3)^3$:

$s^3 + (20K_d-6)s^2 + (20K_p+10)s + 20K_i = s^3 + 9s^2 + 27s + 27$

Equating coefficients:

$20K_d-6 = 9 \Rightarrow K_d = \frac{15}{20} = 0.75$
$20K_p+10 = 27 \Rightarrow K_p = \frac{17}{20} = 0.85$
$20K_i = 27 \Rightarrow K_i = \frac{27}{20} = 1.35$

(c) The impact of implementing the derivative term as $K_d \dot{e}(t)$ versus $-K_d \dot{y}(t)$:\ldots
In summary, using $-K_d \dot{y}(t)$ instead of $K_d \dot{e}(t)$ can provide a more practical implementation of the derivative term in the PID controller, but it may require additional filtering to mitigate high-frequency noise.
    \end{claude3output}

    \begin{humaninput}
    \textbf{Human Input}: Carefully check your solution for Part (c).
    \end{humaninput}
    \begin{claude3output}
    \textbf{Claude 3 Opus}: I apologize for my oversight. Let me take another careful look at Part (c) \ldots
$K_d \dot{e}(t)$ introduces zeros, which can lead to a faster response but may cause overshoot or undershoot.
$-K_d \dot{y}(t)$ affects only the poles, which can improve stability and dampen the response but may result in a slower response.
The choice between the two implementations depends on the specific requirements and constraints of the control system, such as the desired response characteristics, robustness to noise, and the presence of high-frequency components in the reference signal or measurements.
I apologize for my previous oversights and thank you for encouraging me to revisit this problem more carefully.
    \end{claude3output}
\end{conversationbox}

\noindent
For the above problem, both GPT-4 and Gemini 1.0 Ultra cannot give any meaningful values for the PID control gains (Part (b)). In contrast, Claude 3 Opus can solve (a) and (b) in a zero-shot setting. After self-checking, Claude 3 Opus can also give a more precise answer for Part (c). This example again demonstrates the self-correction capabilities of Claude 3 Opus. There are cases where GPT-4 and Gemini 1.0 Ultra also benefit from self-correction, although the resultant improvements  are less significant for these two LLMs.

Obviously, self-correction does not fix all the mistakes made by LLMs.
For instance,  we have observed that self-correction does not work well for problems involving visual elements such as Bode plots and Nyquist plots.
More study on how to combine the self-correction capabilities of LLMs with advanced reasoning search is needed in the future.

\subsection{Sensitivity to the Problem Statements} 
\label{sec:sen}

It is known that LLM responses can be sensitive to the small changes in the human input \cite{sclar2023quantifying,salinas2024butterfly,zou2023universal,guo2024cold}.
Now we discuss some similar phenomenon of LLMs on ControlBench. Specifically, we observed that for some ControlBench problems, very small changes in the problem statements can completely change the LLM responses. We will present two examples to demonstrate this phenomenon. First, we revisit the second example in Section \ref{sec:success} (Problem 5.2 from ControlBench), and show that an almost negligible change in the problem statement can cause Gemini 1.0 Ultra to completely change the solution, as demonstrated below. 
\begin{conversationbox}
\centering{{\bf A Revisit of Problem 5.2 from ControlBench}}
    \begin{humaninput}
    \textbf{Human Input (Original)}: A plant with nominal model 
\begin{equation}\label{eq:Go}
    G_o(s) = \frac{1}{(s + 1)^2}
\end{equation}
is in a feedback loop under control with a PI controller having transfer function $C(s) = \frac{K_ps + K_i}{s}$.
Determine whether this controller configuration can be used to achieve full pole placement.
    \end{humaninput}
\begin{geminioutput}
    \textbf{Gemini 1.0 Ultra}: Here's a breakdown of how to determine if full pole placement is achievable  \ldots Because the number of adjustable parameters from the controller matches the number of poles in the plant, you can likely achieve full pole placement with the given PI controller configuration. Note that this assumes the system is controllable; you'd still need to verify its controllability.
\end{geminioutput}
\begin{humaninput}
    \textbf{Human Input (Modified)}: A plant with nominal model 
$G_o(s) = \frac{1}{(s + 1)^2}$
is in a feedback loop under control with a PI controller having transfer function $C(s) = \frac{K_ps + K_i}{s}$.
Determine whether this controller configuration can be used to achieve full pole placement.
    \end{humaninput}
\begin{geminioutput}
    \textbf{Gemini 1.0 Ultra}: Here's how to determine if full pole placement can be achieved with the given controller configuration  \ldots The given PI controller configuration cannot be used to achieve full pole placement for the plant with the given nominal model \ldots
\end{geminioutput}
\end{conversationbox}

\noindent In the above example, the only change made in the problem statement is that we absorb Equation \eqref{eq:Go} (which originally takes one row) into the text. This is a very minor formatting change. However, the response from Gemini 1.0 Ultra is completely changed. For the modified problem, Gemini 1.0 Ultra gives the correct conclusion. For this specific problem, GPT-4 and Claude 3 Opus are more robust against the changes in the problem statement. Next, we revisit Part (a) of Problem 5.11 from ControlBench (the third example in Section \ref{sec:success}), and demonstrate that our leading LLM Claude 3 Opus sometimes faces a similar issue. 

\begin{conversationbox}
\centering{{\bf Revisiting Part (a) of Problem 5.11 from ControlBench}}
    \begin{humaninput}
    \textbf{Human Input (Original)}: Consider the plant with the following transfer function:
\begin{align*}
    G(s) = \frac{505}{s^3+21s^2+121s+101}
\end{align*}
The poles of this system are at \(s = -1, -10 \pm j\), and the DC gain is $5$.
\begin{itemize}
    \item[(a)] What is the dominant pole approximation \(G_a(s)\) for this plant?
\end{itemize}
    \end{humaninput}
  
\begin{claude3output}
\textbf{Claude 3 Opus}: 
(a) The dominant pole approximation $G_a(s)$ is found by considering the pole closest to the origin (the dominant pole) and neglecting the other poles. In this case, the dominant pole is at $s = -1$. The transfer function of the dominant pole approximation is:
$G_a(s) = \frac{K}{s+1}$,
where $K$ is the DC gain of the original system. Given that the DC gain is 5, we have:
$G_a(s) = \frac{5}{s+1}$.
\end{claude3output}

\begin{humaninput}
    \textbf{Human Input (Modification 1)}: Consider the plant with the following transfer function:
\begin{align*}
    G(s) = \frac{505}{s^3+21s^2+121s+101}
\end{align*}
The poles of this system are at \(s = -1, -10 \pm j\).
\begin{itemize}
    \item[(a)] What is the dominant pole approximation \(G_a(s)\) for this plant?
\end{itemize}
    \end{humaninput}
  
\begin{claude3output}
\textbf{Claude 3 Opus}: 
(a) To find the dominant pole approximation, we need to identify the pole closest to the origin (i.e., the dominant pole) and factor it out from the transfer function \ldots
Now, we can approximate the transfer function by considering only the dominant pole $G_a(s) = \frac{K}{s+1}$\ldots
We can use the residue formula: $\lim_{s\to -1} \frac{505}{s^2+20s+101} 
= \frac{505}{(-1)^2+20(-1)+101} 
= \frac{505}{82} 
\approx 6.159$.
Therefore, the dominant pole approximation is $G_a(s) = \frac{6.159}{s+1}$.
\end{claude3output}
\begin{humaninput}
    \textbf{Human Input (Modification 2)}: Consider the plant with the following transfer function:
\begin{align*}
    G(s) = \frac{505}{s^3+21s^2+121s+101}
\end{align*}
\begin{itemize}
    \item[(a)] Use the DC gain to compute the dominant pole approximation \(G_a(s)\) for this plant.
\end{itemize}
    \end{humaninput}
  
\begin{claude3output}
\textbf{Claude 3 Opus}: 
(a) To find the dominant pole approximation, we first need to find the DC gain and the dominant pole of the system \ldots
Now, we can write the dominant pole approximation as: $G_a(s) = \frac{K}{s-p} = \frac{5}{s+1}$,
where \(K\) is the DC gain and \(p\) is the dominant pole.
\end{claude3output}
\end{conversationbox}

\noindent We tried two modifications for this problem. First, we remove the information of the DC gain being $5$ from the problem statement. Then all of a sudden, Claude 3 Opus fails in recognizing that the DC gain is needed for the dominant pole approximation and makes a weird mistake in giving the wrong answer $G_a(s) = \frac{6.159}{s+1}$. Then one may wonder whether we have to give the exact DC gain information to Claude 3 Opus to ensure that a correct solution can be found. In our second modified problem statement, we remove both the pole and DC gain information, and specifically prompt Claude 3 Opus to use the DC gain for its calculation of the dominant pole approximation. Interestingly, Claude 3 Opus is able to calculate the DC gain and the dominant pole by itself, reaching to the right final answer under this special prompt on ``using the DC gain." 
However, we emphasize that adding the prompt ``DC gain" does not guarantee Claude 3 Opus to solve dominant pole approximation problems.\footnote{For instance, Claude 3 Opus can fail on Problem 1.21 (another problem on dominate pole approximation) from ControlBench, with or without the prompt ``DC gain."}
Overall, the above discussion highlights the need of future research efforts in addressing the consistency and robustness of LLM responses for answering control-related questions.

\section{ControlBench-C: Facilitating Evaluations by Non-Control Experts}
\label{sec:controlbenchC}

The evaluations of LLMs on ControlBench are conducted by a panel of human experts. So far we do not have a fast automated way to evaluate LLM responses on ControlBench. To partially address this issue, we convert 100 problems from ControlBench into single-answer multiple-choice questions, leading to the ControlBench-C dataset. The reason ControlBench-C has fewer problems than ControlBench is that some of the problems from ControlBench involve complicated reasoning and rigorous mathematical proofs, and  it is quite difficult to convert those problems into a multiple-choice format. ControlBench-C  is designed to support evaluations by researchers from diverse backgrounds, including those not specialized in control theory, and has been posted on our project website. Currently, the problems in ControlBench-C are provided by LaTex descriptions, that can be easily converted into the \texttt{Json} format for automated evaluations using API calls. 
For illustrative purposes, consider the following example from ControlBench-C. 

\begin{conversationbox}
\centering{{\bf Modifying Problem 3.5 from ControlBench as a Multi-Choice Problem}}
    \begin{humaninput}
    \textbf{Human Input}: What is the unit step response of a continuous system whose transfer function has a zero at \(-1\), a pole at \(-2\), and a gain factor of \(2\)?
Select the correct option from the following choices:
\begin{enumerate}
    \item [(a)] \(y(t) = 1 + e^{-t}, \quad t \geq 0 \)
    \item [(b)] \(y(t) = 2 + e^{-2t}, \quad t \geq 0 \)
    \item [(c)] \(y(t) = 1 + e^{2t}, \quad t \geq 0 \)
    \item [(d)] \(y(t) = 1 + e^{-2t}, \quad t \geq 0 \)
\end{enumerate}
Please output your choice by just choosing *ONLY ONE* option from (a), (b), (c), or (d), and provide a short explanation below in JSON format by filling in the placeholders in []. Your answer should not include any further explanation, and remember to use parentheses; for example, "(a)" is correct, not "a": 

{

"Choice": "[ONLY ONE choice from (a), (b), (c), (d)]",

"Reason": "[Your explanation]"

}
    \end{humaninput}
\end{conversationbox}

\noindent The above problem is much simpler than the original problem which asks one to compute $y(t)$. For the above problem, all three LLMs choose the correct answer (d). Gemini 1.0 Ultra provides the reason that the system is stable (pole in the left-half plane), and the final value matches the DC gain of the transfer function, while Claude 3 Opus uses the reason that the inverse Laplace transform of the partial fraction expansion gives the unit step response as $y(t) = 1 + e^{-2t}$ for $t \ge 0$.
Interestingly, as discussed in Section~\ref{sec:FM}, Claude 3 Opus uses a similar reasoning path for the original problem (Problem 3.5 from ControlBench) but makes a mistake in the calculations. We can see that somehow such a mistake is avoided when the problem is formulated as a single-answer multi-choice problem.

\begin{table*}[t]
  \centering
  \caption{Accuracy (ACC) and Self-Checked Accuracy (ACC-s) of GPT-4, Claude 3 Opus, and Gemini 1.0 Ultra on ControlBench-C. Both ACC and ACC-s are verified by examining the choice, without considering whether the reasoning is correct or not.  The best
results for each metric are highlighted in bold.}
  \label{tab: main result1}
  \resizebox{1\columnwidth}{!}{%
  \begin{tabular}{l|cc|cc|cc}
    \toprule
     & \multicolumn{2}{c|}{\textbf{GPT-4}} & \multicolumn{2}{c|}{\textbf{Claude 3 Opus}} & \multicolumn{2}{c}{\textbf{Gemini 1.0 Ultra}} \\
   \midrule
   \textbf{Topics} & ACC $\uparrow$ & ACC-s $\uparrow$ & ACC & ACC-s & ACC & ACC-s \\
\midrule
Background & \textbf{73.3\%  (11/15)}  & \textbf{100\% (15/15)}  & 66.7\% (10/15) &80.0\% (12/15)  & \textbf{73.3\% (11/15)}  &80.0\% (12/15) \\
Stability & \textbf{83.3\% (10/12)}  & \textbf{91.7\% (11/12)}  & 50.0\% (6/12)  &\textbf{91.7\% (11/12)}  & \textbf{83.3\% (10/12)}  & \textbf{91.7\% (11/12)}  \\
Time response & 76.4\% (13/17) &76.4\% (13/17) & \textbf{82.3\% (14/17)}  &\textbf{94.1\% (16/17)}  & 64.7\% (11/17) & \textbf{94.1\% (16/17)}  \\
Block diagrams &\textbf{50.0\% (1/2)} &\textbf{50.0\% (1/2)} &\textbf{50.0\% (1/2)}  &\textbf{50.0\% (1/2)}  & 0.0\% (0/2)  & \textbf{50.0\% (1/2)}  \\
Control System Design & \textbf{43.7\% (7/16)} &50.0\% (8/16) & 37.5\% (6/16)  &\textbf{56.2\% (9/16)}  & \textbf{43.7\% (7/16)}  & \textbf{56.2\% (9/16)} \\
Bode Analysis &36.3\% (4/11) & \textbf{90.9\% (10/11)} & \textbf{54.5\% (6/11)}  &\textbf{90.9\% (10/11)}  & 36.3\% (4/11)  &63.6\% (7/11)\\
Root-Locus Design &40.0\% (2/5) &40.0\% (2/5) & \textbf{80.0\% (4/5)} &\textbf{100\% (5/5)}  & 60.0\% (3/5)  &60.0\% (3/5) \\
Nyquist Design & 25.0\% (1/4) &50.0\% (2/4) &\textbf{50.0\% (2/4)}  &\textbf{75.0\% (3/4)} & 0.0\% (0/4)  &50.0\%  (2/4)\\
Gain/Phase Margins &\textbf{71.4\% (5/7)} &\textbf{85.7\% (6/7)} & 57.1\% (4/7)  &\textbf{85.7\% (6/7)}  & 57.1\% (4/7)  & \textbf{85.7\% (6/7)} \\
System Sensitivity Measures  &\textbf{100.0\% (3/3)} &\textbf{100.0\% (3/3)} & \textbf{100.0\% (3/3)}  &\textbf{100.0\% (3/3)}  & 66.7\% (2/3)  &\textbf{100.0\% (3/3)} \\
Loop-shaping  &0.0\% (0/1) &0.0\% (0/1) & 0.0\% (0/1)  &\textbf{100\% (1/1)}  & \textbf{100\% (1/1)}  & \textbf{100\% (1/1)} \\
Advanced Topics & \textbf{100\% (7/7)} & \textbf{100\% (7/7)} & 42.9\% (3/7)  &85.7\% (6/7)  & 42.9\% (3/7)  &57.1\% (4/7) \\
\midrule 
Total & \textbf{64.0\% (64/100)} &78.0\% (78/100) & 59.0\% (59/100)  &\textbf{83.0\% (83/100)}  & 56.0\% (56/100)  &75.0\% (75/100)  \\ 
\bottomrule
\end{tabular}
}
\end{table*}

For completeness, we also evaluate the ACC and ACC-s of GPT-4, Claude 3 Opus, and Gemini 1.0 Ultra on ControlBench-C. For fast automated evaluations, both ACC and ACC-s are calculated via examining the choices from the LLMs, without considering whether the reasoning is correct or not. Such evaluations can be easily done by non-control experts. However, the downside is that sometimes LLMs can pick the right choice based on wrong reasoning, causing concerns about whether such testing reflects the true capabilities of LLMs in solving control problems. For ControlBench-C, GPT-4 and Claude 3 Opus have similar performances. GPT-4 achieves the highest ACC, while Claude 3 Opus achieves the highest ACC-s. 
We view ControlBench-C as a complement to ControlBench. Specifically, ControlBench-C enables automated evaluations even by non-control experts. However, the results from ControlBench-C do not provide the same level of insight as our previous analysis on ControlBench. An important future task is to develop automated evaluation methods for ControlBench. 

\section{Conclusion and Future Work}
In this paper, we study the capabilities of large language models (LLMs) including GPT-4, Claude 3 Opus, and Gemini 1.0 Ultra in solving undergraduate control engineering problems. To support the study, we introduce a benchmark dataset, {\bf ControlBench}. We offer comprehensive insights from control experts to uncover the current potential and limitations of LLMs. We believe that our work is just a starting point for further studies of LLM-based methods for control engineering.
 We conclude our paper with a brief discussion on future research directions.

\paragraph{Expansion of the problem set.} Increasing the diversity of the problem set in the  ControlBench dataset could provide a more comprehensive evaluation of LLMs.  By introducing more complex and challenging control design tasks, the dataset can push the capabilities of LLMs further. It will also be interesting to include more problems on advanced topics such as nonlinear control \cite{khalil01}, stochastic control \cite{aastrom2012introduction}, robust control \cite{zhou1996robust}, adaptive control \cite{aastrom2013}, and quantum control \cite{d2021introduction}.

\paragraph{Control-oriented prompting.} It is well known that  LLMs can often be steered to generate prescribed responses as long as one finds the right prompts
\cite{bhargava2023s,arora2022ask,geiping2024coercing}. 
In many situations, short prompts, referred to as ``magic words" \cite{bhargava2023s}, are sufficient for the purpose of steering LLMs. It will be interesting to investigate whether one can combine domain knowledge and automatic prompt search methods to develop control-oriented prompts and instructions for improving the capabilities of LLMs in control engineering.

\paragraph{Improving reasoning capabilities and tool use abilities for consistency and accuracy.} 
In the future, we also plan to adopt advanced planning strategies such as trees-of-thought prompting \cite{yao2024tree} or Monte Carlo Tree Search (MCTS) \cite{hao2023reasoning} and integrate external resources like programming code or \texttt{Matlab} control toolboxes \cite{he2023solving, hao2024toolkengpt}. Those strategies aim at enhancing the model accuracy and consistency in reasoning and calculations for the control design task. One main issue observed in our paper is that LLMs can sometimes generate inconsistent answers even for the same problem. Improving the reasoning and calculation capabilities of LLMs can potentially lead to solutions for this issue.

\paragraph{Efficient evaluation.} Efficient evaluation methods are crucial for scaling the application of LLMs in control engineering. Automating the evaluation process, while ensuring it complements expert human assessments, could streamline the validation of LLM outputs and facilitate more rapid advancements in the field. There is a need to develop new methods and metrics for fast automated evaluation of LLMs on control benchmarks.

\paragraph{Vision-language models for handling various plots.}
As observed in our paper, GPT-4, Claude 3 Opus, and Gemini 1.0 Ultra all have difficulties in solving problems involving Bode plots, Nyquist plots, and root locus plots. Therefore, it is important to develop new methods for incorporating control-domain knowledge into state-of-the-art vision-language models \cite{dai2024instructblip, li2023blip, liu2024visual, zhu2023minigpt}.

\section*{Acknowledgement}

U. Syed, X. Guo, A. Havens, and B. Hu are generously supported by the NSF award CAREER-2048168.

\section*{Potential Social Impact}

The aspect of AI safety in the context of integrating large language models (LLMs) into control engineering is paramount, especially given the potential for these models to be applied in critical infrastructure and systems.  As we look towards a future where LLMs may play a significant role in designing, optimizing, and maintaining control systems, we must prioritize the development of safety protocols and standards to govern their deployment.
The integration of LLMs in control engineering also raises important ethical considerations. As these models begin to influence decision-making in control systems, questions regarding accountability, transparency, and the potential for unintended consequences must be addressed. Developing frameworks that clearly delineate the responsibilities of human operators and LLMs will be crucial. Additionally, ensuring that LLMs are designed with fairness and bias mitigation in mind will help prevent the propagation of existing prejudices into control engineering solutions.
To address these challenges and opportunities, fostering a collaborative environment between control engineers, AI researchers, ethicists, and policymakers is essential. Such collaborations can lead to the development of interdisciplinary solutions that not only enhance the technical capabilities of LLMs in control engineering but also ensure that their deployment is safe, ethical, and socially beneficial. By combining domain-specific knowledge with advancements in AI, we can create robust frameworks for the responsible use of LLMs in control engineering.
In addition, the development of comprehensive regulatory frameworks and standards specific to the use of LLMs in control engineering will be crucial. These frameworks should address aspects such as the validation of LLM outputs, the ethical use of AI in engineering applications, and the safety of AI-driven control systems. Establishing clear guidelines and standards will not only promote the safe and responsible use of LLMs but also foster public trust in AI-enhanced control engineering solutions.

The integration of large language models (LLMs) into control engineering, and their broader application across various disciplines, also prompts a critical examination of their potential negative social impacts, particularly in the realm of education. The accessibility and efficiency of LLMs in solving complex problems might inadvertently lead to a reliance on these tools among students, potentially undermining the development of foundational problem-solving skills and critical thinking. This scenario could lead to a superficial understanding of complex subjects, diminishing the educational process's depth and rigor.
To mitigate these potential negative impacts on education, it is crucial to adopt a balanced approach that leverages the benefits of LLMs while fostering deep learning and skill development. One effective strategy could involve integrating LLMs into the curriculum as supplementary tools rather than primary problem solvers. Educators can design assignments and projects that require students to critically evaluate LLM-generated solutions, encouraging deeper engagement with the material and promoting critical thinking.
Moreover, developing educational frameworks that emphasize the understanding of underlying principles rather than solely focusing on obtaining solutions can help maintain the educational quality. Incorporating project-based learning, where students must apply concepts to real-world scenarios, can ensure that they develop practical skills and a comprehensive understanding of the subject matter.
Finally, ethical training regarding the use of LLMs and other AI tools in academic settings may be integrated into curricula to instill a sense of responsibility and integrity among students.

\bibliographystyle{plainnat}
\bibliography{main}

\end{document}